\input{epsf}

\newcounter{letter}     
\newenvironment{alphalist}{\begin{list}
{{\normalshape(\alph{letter})}}{\usecounter{letter}}}{\end{list}}

        \newcommand{\be}{\begin{equation}}
        \newcommand{\ee}{\end{equation}}
        \newcommand{\ba}{\begin{eqnarray}}
        \newcommand{\ea}{\end{eqnarray}}
        \newcommand{\ban}{\begin{eqnarray*}}
        \newcommand{\ean}{\end{eqnarray*}}
        \newcommand{\barr}{\begin{array}}
        \newcommand{\earr}{\end{array}}

\def\=>{\Rightarrow}
\newcommand{\To}{\Rightarrow}
\def\mapright#1{\smash{\mathop{\longrightarrow}\limits^{#1}}}

\newcommand{\C}{{\Bbb C}}
\newcommand{\et}{\hspace{-0.08in}{\bf .}\hspace{0.1in}}
\newcommand{\BOX}{\hbox {$\sqcap$ \kern -1em $\sqcup$}}

\renewcommand{\hom}{{\rm hom}}

\renewcommand{\to}{\rightarrow}

\newcommand{\maps}{\colon}

\newcommand{\iso}{\cong}

\newcommand{\elt}{{\rm elt}}

\newtheorem{thm}{Theorem}    

\newtheorem{defn}[thm]{Definition}

	\documentstyle[12pt,amsfonts,diagram]{article}
\begin{document}

\hfuzz=9pt      

      \begin{center}
      {\bf An Introduction to $n$-Categories\\}
       \vspace{0.3cm}
      {\em John C.\ Baez \\}
      \vspace{0.3cm}
      {\small Department of Mathematics,  University of California\\ 
      Riverside, California 92521 \\
      USA\\ }
      \vspace{0.3cm}
      {\small email: baez@math.ucr.edu\\}
      \vspace{0.3cm}
      {\small May 13, 1997 \\ }
      \end{center}

\begin{abstract} An $n$-category is some sort of
algebraic structure consisting of objects, morphisms between objects,
2-morphisms between morphisms, and so on up to $n$-morphisms, together
with various ways of composing them.  We survey various concepts of
$n$-category, with an emphasis on `weak' $n$-categories, in which all
rules governing the composition of $j$-morphisms hold only up to
equivalence.  (An $n$-morphism is an equivalence if it is
invertible, while a $j$-morphism for $j < n$ is an equivalence if it is
invertible up to a $(j+1)$-morphism that is an equivalence.)  We discuss
applications of weak $n$-categories to various subjects including
homotopy theory and topological quantum field theory, and review the
definition of weak $n$-categories recently proposed by Dolan and the
author.  \end{abstract}

\section{Introduction}

Very roughly, an $n$-category is algebraic structure consisting
of a collection of `objects', a collection of `morphisms' between
objects, a collection of `2-morphisms' between morphisms, and so on up
to $n$, with various reasonable ways of composing these $j$-morphisms.
A $0$-category is just a set, while a $1$-category is just a category.
Recently $n$-categories for arbitrarily large $n$ have begun to play an
increasingly important role in many subjects.  The reason is that they
let us {\it avoid mistaking isomorphism for equality}.  

In a mere set, elements are either the same or different; there is no
more to be said.  In a category, objects can be different but still `the
same in a way'.  In other words, they can be unequal but still
isomorphic.  Even better, we can explicitly keep track of the way they
are the same: the isomorphism itself.  This more nuanced treatment of
`sameness' is crucial to much of mathematics, physics, and computer
science.  For example, it underlies the modern concept of symmetry:
since an object can be `the same as itself in different ways', it has a
symmetry group, its group of automorphisms.  Unfortunately, in a
category this careful distinction between equality and isomorphism
breaks down when we study the morphisms.  Morphisms in a category are
either the same or different; there is no concept of isomorphic
morphisms.  In a 2-category this is remedied by introducing 2-morphisms
between morphisms.  Unfortunately, in a 2-category we cannot speak of
isomorphic 2-morphisms.  To remedy this we need the notion of
3-category, and so on.

The plan of this paper is as follows.  We do not begin by defining
$n$-categories.  Many definitions have been proposed.  So far, all of
them are a bit complicated.  Ultimately a number of them should turn out
to be equivalent, but this has not been shown yet.  In fact, the correct
sense of `equivalence' here is a rather subtle issue, intimately linked
with $n$-category theory itself.  Thus before mastering the details of
any particular definition, it is important to have a sense of the issues
involved.  Section 2 starts with a rough sketch of various approaches to
defining $n$-categories.  Section 3 describes how $n$-categories are
becoming important in a variety of fields, and how they should allow us
to formalize some previously rather mysterious analogies between
different subjects.  Section 4 sketches a particular definition of `weak
$n$-categories' due to Dolan and the author \cite{BD2}.  In the Conclusions
we discuss the sense in which various proposed definitions of 
weak $n$-category should be equivalent.

\section{Various Concepts of $n$-Category}

To start thinking about $n$-categories it is helpful to use pictures. 
We visualize the objects as 0-dimensional, i.e., points.  We
visualize the morphisms as 1-dimensional, i.e., intervals, or
more precisely, arrows going from one point to another.  In this picture,
composition of morphisms corresponds to gluing together an arrow
$f \maps x \to y$ and an arrow $g \maps y \to z$ to obtain an arrow
$fg \maps x \to z$:

\medskip
\centerline{\epsfysize=0.3in\epsfbox{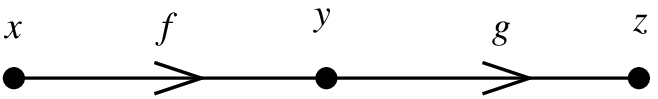}}
\medskip

\noindent Note that while the notation $fg$ for the composite of 
$f \maps x \to y$ and $g \maps y \to z$ is somewhat nonstandard, it 
fits the picture better than the usual notation.

Continuing on in this spirit, we visualize the 2-morphisms as
2-dimensional, and compose 2-morphisms in a way that corresponds
to gluing together 2-dimensional
shapes.  Of course, we should choose some particular shapes for our
2-morphisms.  For example, we could use a `bigon':

\medskip
\centerline{\epsfysize=1.0in\epsfbox{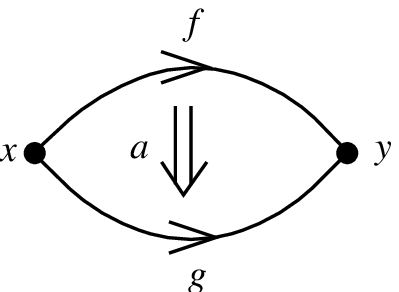}}
\medskip

\noindent as the shape of a 2-morphism $a \maps f \To g$ between
morphisms $f,g \maps x \to y$ with the same source and target.  This is
the sort of 2-morphism used in the standard definitions of `strict
2-categories' \cite{KS} --- usually just called 2-categories --- and the
somewhat more general `bicategories' \cite{Benabou}.  There are two
geometrically natural ways to compose 2-morphisms shaped like bigons.
First, given 2-morphisms $a \maps f \To g$ and $b \maps g \To h$ as
below, we can `vertically' compose them to obtain a 2-morphism $a \cdot
b\maps f \To g$:

\medskip
\centerline{\epsfysize=1.5in\epsfbox{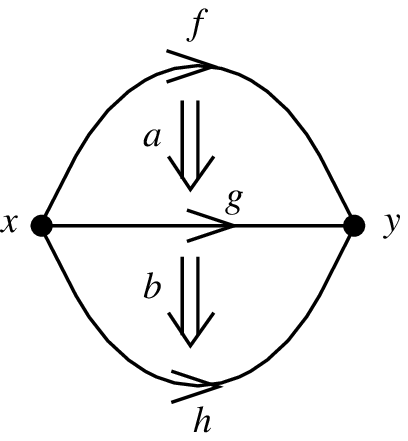}}
\medskip

\noindent Second, given 2-morphisms $a \maps f \To g$ and $b \maps h \To
i$ as below, we can `horizontally' compose them to obtain a 2-morphism
$ab \maps fh \To gi$:

\medskip
\centerline{\epsfysize=1.0in\epsfbox{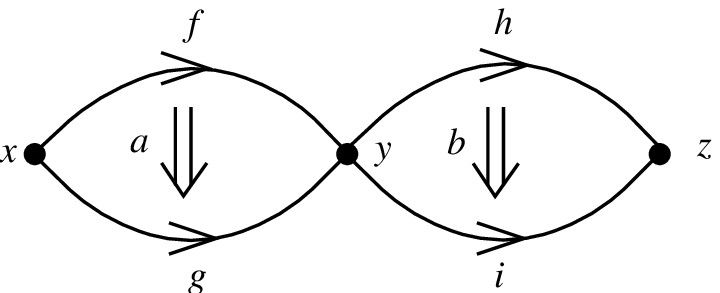}}
\medskip

The definition of a strict 2-category is easy to state.  The objects and
morphisms must satisfy the usual rules holding in a category, while
horizontal and vertical composition satisfy some additional axioms:
vertical and horizontal composition are associative, and for each
morphism $f$ there is a 2-morphism $1_f\maps f \To f$ that is an
identity for both vertical and horizontal composition.  Finally, we
require the following `interchange law' relating vertical and horizontal
composition:
\[    
(a \cdot b)(c \cdot d) = (ac)\cdot (bd)
\]
whenever either side is well-defined.
This makes the following composite 2-morphism unambiguous:

\medskip
\centerline{\epsfysize=1.0in\epsfbox{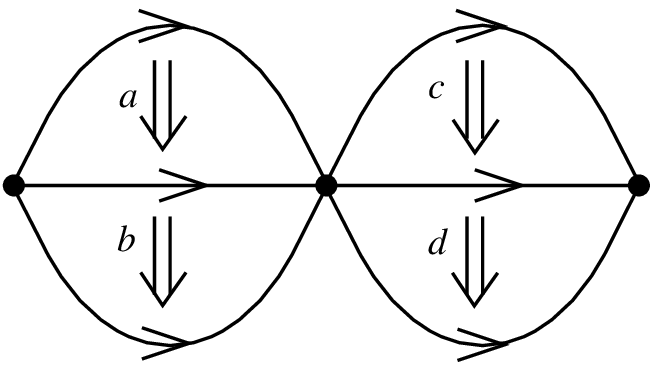}}
\medskip

\noindent We can think of it either as the result of first doing two
vertical composites and then one horizontal composite, or as the result
of first doing two horizontal composites and then one vertical
composite.

The definition of a bicategory is similar, but instead of requiring that
the associativity and unit laws for morphisms hold `on the nose'
as equations, one requires merely that they hold up to isomorphism.
Thus one has invertible `associator' 2-morphisms 
\[        a_{f,g,h} \maps (fg)h \To f(gh) \]
for every composable triple of morphisms, as well as invertible
2-morphisms called `left and right identity constraints'
\[       l_f \maps 1_x f \To f, \qquad
          r_f \maps f 1_y \To f \] 
for every morphism $f \maps x \to y$.  These must satisfy some equations
of their own.  For example, repeated use of associator lets one go from
any parenthesization of a product of morphisms to any other
parenthesization, but one can do so in many ways.  To ensure that all
these ways are equal, one imposes the Stasheff `pentagon identity',
which says that the following diagram commutes: 
\[
\begin{diagram}[(fg)(hi)] 
\node{((fg)h)i} \arrow{e,t}{a_{fg,h,i}}\arrow{s,t}{a_{f,g,h}1_i} 
\node{(fg)(hi)} \arrow{e,t}{a_{f,g,hi}}
\node{f(g(hi))} \\ 
\node{(f(gh))i} \arrow[2]{e,t}{a_{f,gh,i}}
\node[2]{f((gh)i)} \arrow{n,r}{1_f a_{g,h,i}} \end{diagram} 
\]
Mac Lane's coherence theorem \cite{Maclane} says that this identity
suffices.  Similarly, given morphisms $f \maps x \to y$ and $g \maps y
\to z$, one requires that the following triangle commute:
\[  
\begin{diagram}[f(1_yg)]
\node{(f1_y)g} \arrow[2]{e,t}{a_{f,1_y,g}} \arrow{se,b}{r_fg} 
\node[2]{f(1_yg)} \arrow{sw,r}{fl_g} \\
\node[2]{fg} 
\end{diagram}
\]
One also requires that the associators and unit constraints are natural
with respect to their arguments.  Also, as with strict 2-categories,
one requires that vertical composition be associative, that vertical and
horizontal composition satisfy the interchange law, and that the 
morphisms $1_x$ are identities for vertical composition.

While bicategories at first seem more clumsy than strict 2-categories,
they arise more often in applications.  The reason is that in
applications typically `everything is only true up to something'.  In
a sense, the whole point of introducing $(n+1)$-morphisms is to allow
$n$-morphisms to be isomorphic rather than merely equal.  From this
point of view, it was inappropriate to have imposed equational laws
between 1-morphisms in the definition of a strict 2-category, and the
definition of bicategory corrects this problem.  This is known as
`weakening'.  

To see some bicategories that are not strict 2-categories, consider
bicategories with one object.  Given a bicategory $C$ with one object
$x$, we can form a category $\tilde C$ whose objects are the morphisms
of $C$ and whose morphisms are the 2-morphisms of $C$.  This is a
special sort of category: we can `multiply' the objects of $\tilde C$,
since they are really just morphisms in $C$ from $x$ to itself.  We call
this sort of category --- one that is really just a bicategory with one
object --- a `weak monoidal category'.  We can do the same thing
starting with a strict 2-category and get a `strict monoidal category'.

The category Set becomes a weak monoidal category if we multiply sets 
using the Cartesian product.  However, it is not a strict monoidal
category!  The reason is that the Cartesian product is not strictly
associative:
\[               (X \times Y) \times Z \ne X \times (Y \times Z) .\]
To see this, one needs to pry into the set-theoretic definition of
ordered pairs.  The usual von Neumann definition is $(x,y) = \{\{x\},
\{x,y\}\}$, and using this, we clearly do not have strict associativity
for the Cartesian product.  Instead, we have associativity {\it up to
a specified isomorphism}, the associator:
\[       a_{X,Y,Z} \maps (X \times Y) \times Z \to X \times (Y \times Z)
\]
which satisfies the pentagon identity.  This is a typical example of
how the bicategories found `in nature' tend not to be strict 2-categories. 

So far we have only considered bigons as possible shapes for
2-morphisms, but there are many other choices.  For example, we might
wish to use triangles going from a pair of morphisms $f \maps x \to y$,
$g \maps y \to z$ to a morphism $h \maps x \to z$:

\medskip
\centerline{\epsfysize=1.2in\epsfbox{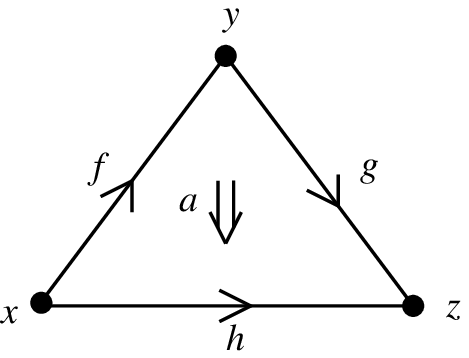}}
\medskip

\noindent  There is no good way to glue together triangles of
this type to form other triangles of this type, but if we also 
allow the `reverse' sort of triangle going from a single morphism 
to a pair: 

\medskip
\centerline{\epsfysize=1.2in\epsfbox{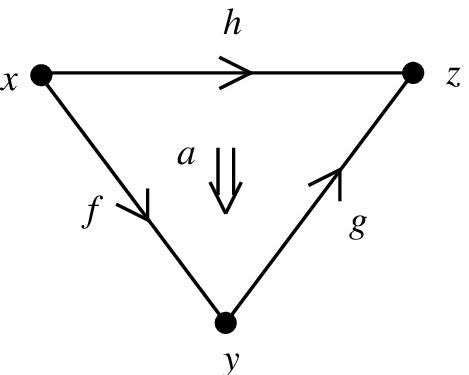}}
\medskip

\noindent then there are various ways to glue together 3 triangles to
form a larger one.  For example:

\medskip
\centerline{\epsfysize=1.2in\epsfbox{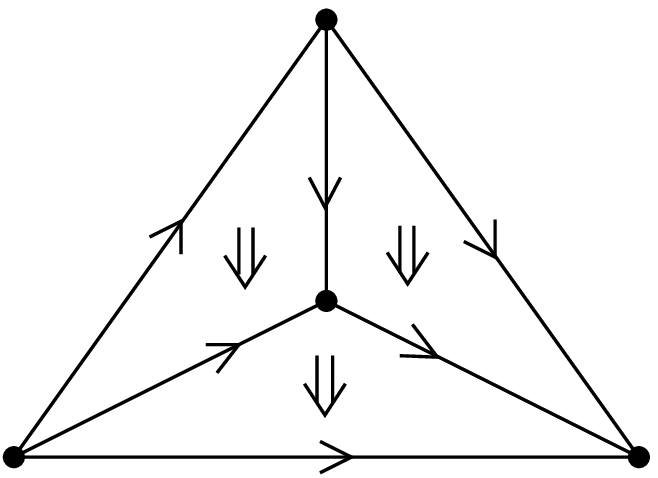}}
\medskip

\noindent One can do quite a bit of topology in a purely combinatorial
way using triangles and their higher-dimensional analogues, called
simplices \cite{May}.  In these applications one often assumes that all
$j$-morphisms are invertible, at least in some weakened sense.  In the
2-dimensional case, this motivates the idea of `reversing' a triangle.

Another approach would be to use squares going from a pair of
morphisms $f \maps x \to y$, $g \maps y \to z$ to a pair $h \maps x \to
w$, $i \maps w \to z$:

\medskip
\centerline{\epsfysize=1.2in\epsfbox{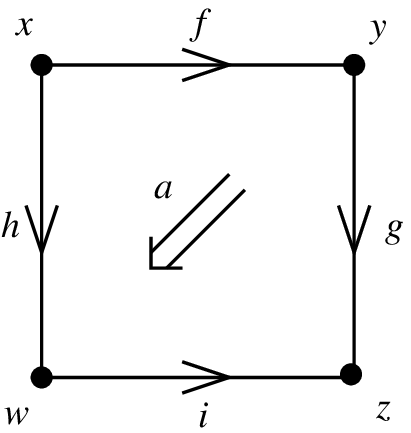}}
\medskip

Much as with bigons, one can compose squares vertically and horizontally
and require an `interchange law' relating these two types of
composition.  This is the idea behind the definition of `double
categories' \cite{Ehresmann,KS}, where the vertical arrows are treated
as of a different type as the horizontal ones.  If one treats the
vertical and horizontal arrows as the same type, one obtains a theory
equivalent to that of strict 2-categories.

Alternatively, one might argue that the business of picking a particular
shape of 2-morphism as `basic' is somewhat artificial.  One might
instead allow all possible polygons as shapes for 2-morphisms.
The idea would be to use polygons whose boundary is divided into two parts
having the arrows consistently oriented:

\vbox{
\medskip
\centerline{\epsfysize=1.2in\epsfbox{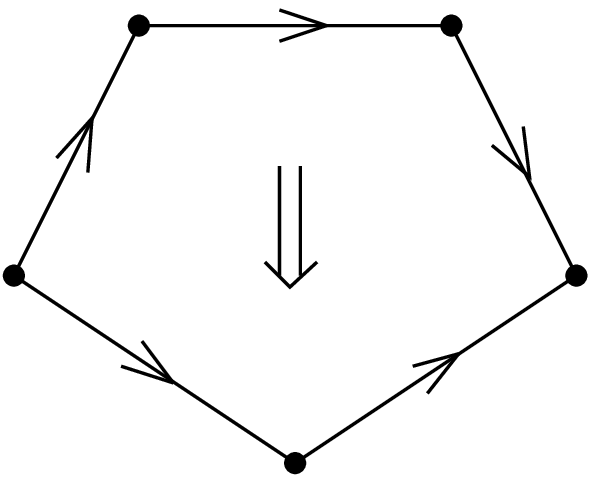}}
\medskip
}

\noindent There are many ways to compose such polygons.  However, while this
approach might seem more general, one can actually define and work
with these more general polygons within the theory of 
strict 2-categories \cite{Johnson,Power,Power2}.  

Yet another approach would be to use only polygons having many `inface'
but only one `outface', like this:

\medskip
\centerline{\epsfysize=1.0in\epsfbox{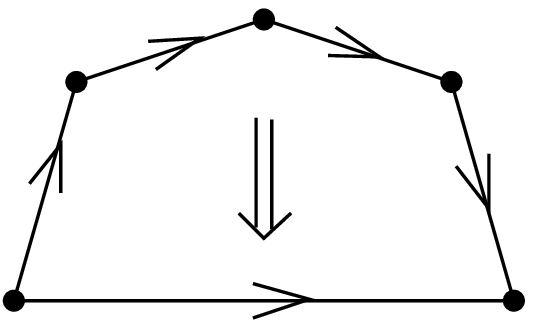}}
\medskip

\noindent As we shall see, this has certain advantages of
its own.

What about $n$-categories for higher $n$?  In general, $j$-morphisms can
be visualized as $j$-dimensional solids, and part of their boundary
represents the source while the rest represents the target.  However,
as $n$ increases one can imagine more and more definitions of
`$n$-category', because there there are more and more choices for the
shapes of the $j$-morphisms.  The higher-dimensional analogues of bigons
are called `globes'.  Globes are the basic shape in the traditional
approach to `strict $n$-categories' \cite{EK,Crans}.  There has also
been a lot of work on globular weak $n$-categories, such as Gordon,
Power and Street's `tricategories' \cite{GPS}, Trimble's
`tetracategories' \cite{Trimble}, and Batanin's `weak
$\omega$-categories' \cite{Batanin}, which can have $j$-morphisms of
arbitrarily high dimension.  The higher-dimensional analogues of
triangles, namely simplices, are used in the `Kan complexes' favored by
topologists \cite{May}, as well as in Street's `simplicial weak
$\omega$-categories' \cite{Street} and Lawrence's `$n$-algebras'
\cite{Lawrence}.  The higher-dimensional analogues of squares, namely
cubes, are used in Ehresmann's `$n$-tuple categories' \cite{Ehresmann},
as well as the work of Brown and his collaborators \cite{Brown}.
Finally, Dolan and the author \cite{BD2} have given a definition of
`weak $n$-categories' based on some new shapes called `opetopes'.  We
describe these in Section 4.

In addition to the issue of shapes for $j$-morphisms, there is the issue
of the laws that composition operations should satisfy.  Most
importantly, there is the distinction between `strict' and `weak'
approaches.  In the `strict' approach, composition of $j$-morphisms
satisfies equational laws for all $j$.  The philosophy behind the `weak'
approach is that equations should hold only at the top level,
between $n$-morphisms.   Laws concerning $j$-morphisms for $j < n$
should always be expressed as $(j+1)$-morphisms, or more precisely, 
`equivalences'.  Roughly, the idea here is that an equivalence between
$(n-1)$-morphisms is an invertible $n$-morphism, while an equivalence
between $j$-morphisms for lesser $j$ is recursively defined as a
$(j+1)$-morphism that is invertible {\it up to equivalence}.  

Strict $n$-categories are fairly well-understood \cite{Crans}, but the
interesting and challenging sort of $n$-categories are the weak ones.
Weak $n$-categories are interesting because these are the ones that tend
to arise naturally in applications.  The reason for this is simple yet
profound.  Equations of the form $x = x$ are completely useless.  All
interesting equations are of the form $x = y$.  Equations of this
form can always be viewed as asserting the existence of a reversible
sort of computation transforming $x$ to $y$.  In $n$-categorical terms,
they assert the existence of an equivalence $f \maps x \to y$.  To face
up to this fact, it is helpful to systematically avoid equational laws
and work explicitly with equivalences, instead.  This leads naturally to
working with weak $n$-categories, and eventually weak $\omega$-categories.

The reason $n$-categories are challenging is that when equational laws
are replaced by equivalences, these equivalences need to satisfy new
laws of their own, called `coherence laws', so one can manipulate them
with some of the same facility as equations.  The main problem of weak
$n$-category theory is: how does one systematically determine these
coherence laws?  A systematic approach is necessary, because in general
these coherence laws must themselves be treated not as equations but as
equivalences, which satisfy further coherence laws of their own, and so
on!  This quickly becomes very bewildering if one proceeds on an ad hoc
basis.

For example, suppose one tries to write down definitions of `globular
weak $n$-categories', that is, weak $n$-categories in the approach where
the $j$-morphisms are shaped like globes.  These are usually called
categories, bicategories, tricategories, tetracategories, and so on.
The definition of a category is quite concise; the most complicated
axiom is the associative law $(fg)h = f(gh)$.  As we
have seen, in the definition of a bicategory this law is replaced by a
2-morphism, the associator, which in turn satisfies the pentagon
identity.  In the definition of a tricategory, the pentagon identity is
replaced by a 3-isomorphism satisfying a coherence law which is best
depicted using a 3-dimensional commutative diagram in the shape of the
3-dimensional `associahedron'.  In the definition of a tetracategory,
this becomes a 3-morphism which satisfies a coherence law given by the
4-dimensional associahedron.  In fact, the associahedra of all
dimensions were worked out by Stasheff \cite{Stasheff} in 1963 using
homotopy theory.  However, there are other sequences of coherence laws
to worry about, spawned by the equational laws of the form $1f = f =
f1$, and also the interchange laws governing the various
higher-dimensional analogues of `vertical' and `horizontal' composition.

At this point the reader can be forgiven for wondering if the rewards of
setting up a theory of weak $n$-categories really justify the labor
involved.  Before proceeding, let us describe some of the things
$n$-categories should be good for.

\section{Applications of $n$-Categories}

One expects $n$-categories to show up in any situation where there are
things, processes taking one thing to another, `meta-processes' taking
one process to another, `meta-meta-processes', and so on.  Clearly
computer science is deeply concerned with such situations.
Some other places where applications are evident include:
\begin{enumerate} 
\item $n$-category theory 
\item homotopy theory 
\item topological quantum field theory 
\end{enumerate} 
The first application is circular, but not viciously so.  The point is
that the study of $n$-categories leads to applications of
$(n+1)$-category theory.  The other two applications may sound abstruse
and specialized, but there is a good reason for discussing them here.
Pure $n$-category theory treats the most general iterated notion of
process.  Homotopy theory limits its attention to processes that are
`invertible', at least up to equivalence.  Topological quantum field
theory focuses attention on processes which have `adjoints' or `duals'.
While generally not invertible even up to equivalence, such processes
are reversible in a broader sense (the classic examples from category
theory being adjoint functors).  In what follows we briefly summarize
all three applications in turn.

\subsection{$n$-Category Theory}

While self-referential, this application is perhaps the most
fundamental.  A 0-category is just a {\it set}.  When one studies sets
one is naturally led to consider the set of all sets.  However, this
turns out to be a bad thing to do, not merely because of Russell's
paradox (which is easily sidestepped), but because one is interested not
just in sets but also in the functions between them.  What is
interesting is thus the {\it category} of all sets, Set.

This category is in some sense the primordial category.  Indeed, the
Yoneda embedding theorem shows how every category can be thought of as a
category of `sets with extra structure'.  However, when we study
categories of sets with extra structure, it turns out to be worthwhile
to develop category theory as a subject in its own right.  In addition
to categories and functors, natural transformations play a crucial role
here.  Thus one is led to study the {\it 2-category} of all categories,
Cat.  This 2-category has categories as objects, functors between
categories as morphisms, and natural transformations between functors as
2-morphisms.

The ladder of $n$-categories continues upwards in this way.  For each
$n$ there is an $(n+1)$-category of all $n$-categories, $n$Cat.  To
really understand $n$-categories we need to understand this
$(n+1)$-category.  Eventually this requires an understanding of
$(n+1)$-categories in general, which then leads us to define $(n+1)$Cat.

There are some curious subtleties worth noting here, though.  The
2-category Cat happens to be a strict 2-category.  We could think of it
as a bicategory if we wanted, but weakening happens not to be needed
here, since functors compose associatively `on the nose', not just up to
a natural transformation.  Using the fact that Cat is the primordial
2-category, one can show that every bicategory is equivalent to a
strict 2-category in a certain precise sense.  Technically speaking, one
proves this using the Yoneda embedding for bicategories \cite{GPS}.

The fact that every weak 2-category can be `strictified' seems to have
held back work on weak $n$-categories: it raised the hope that every
weak $n$-category might be equivalent to a strict one.  It turned out,
however, that the strict and weak approaches diverge as we continue to
ascend the ladder of $n$-categories.  On the one hand, we can always
construct a strict $(n+1)$-category of strict $n$-categories.  On the
other hand, we can construct a weak $(n+1)$-category of weak
$n$-categories.  The latter is {\it not} equivalent to a strict
$(n+1)$-category for $n \ge 2$.

Consider the case $n = 2$.  On the one hand, we can form a strict
3-category 2Cat whose objects are strict 2-categories.  We can
visualize a strict 2-category as a bunch of points, arrows and bigons.
For simplicity, let us consider a very small 2-category $C$ with just
one interesting 2-morphism:

\medskip
\centerline{\epsfysize=1.0in\epsfbox{bigon.eps}}
\medskip

\noindent (We have not drawn the identity morphisms and 2-morphisms.)
The morphisms in 2Cat are called `2-functors'.  A 2-functor $F \maps 
C \to D$ sends objects to objects, morphisms to morphisms, 
and 2-morphisms to 2-morphisms, strictly preserving all structure.  
We can visualize $F$ as creating a picture of the 
2-category $C$ in the 2-category $D$:

\medskip
\centerline{\epsfysize=1.0in\epsfbox{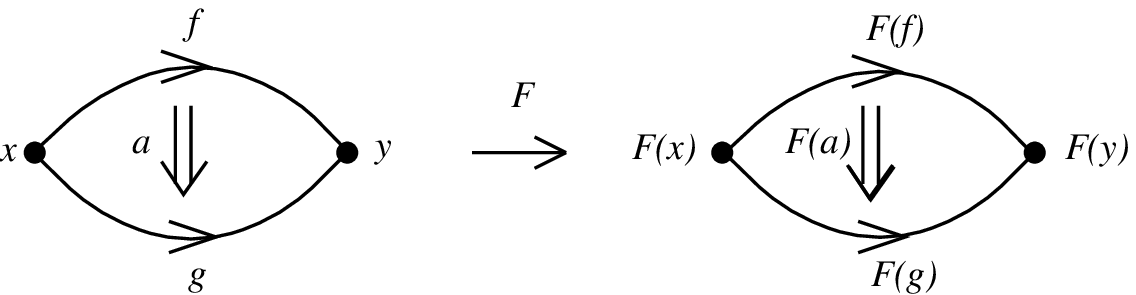}}
\medskip

\noindent The 2-morphisms in 2Cat are called `natural transformations'.  
A natural transformation $A \maps F \To G$ between 2-functors 
$F,G \maps C \to D$ sends each object in $C$ to a morphism in $D$
and each morphism in $C$ to a 2-morphism in $D$, and satisfies
some conditions similar to those in the definition of a natural transformation 
between functors.   We can visualize $A$ as a prism going from one picture
of $C$ in $D$ to another, built using commutative squares:

\medskip
\centerline{\epsfysize=2.0in\epsfbox{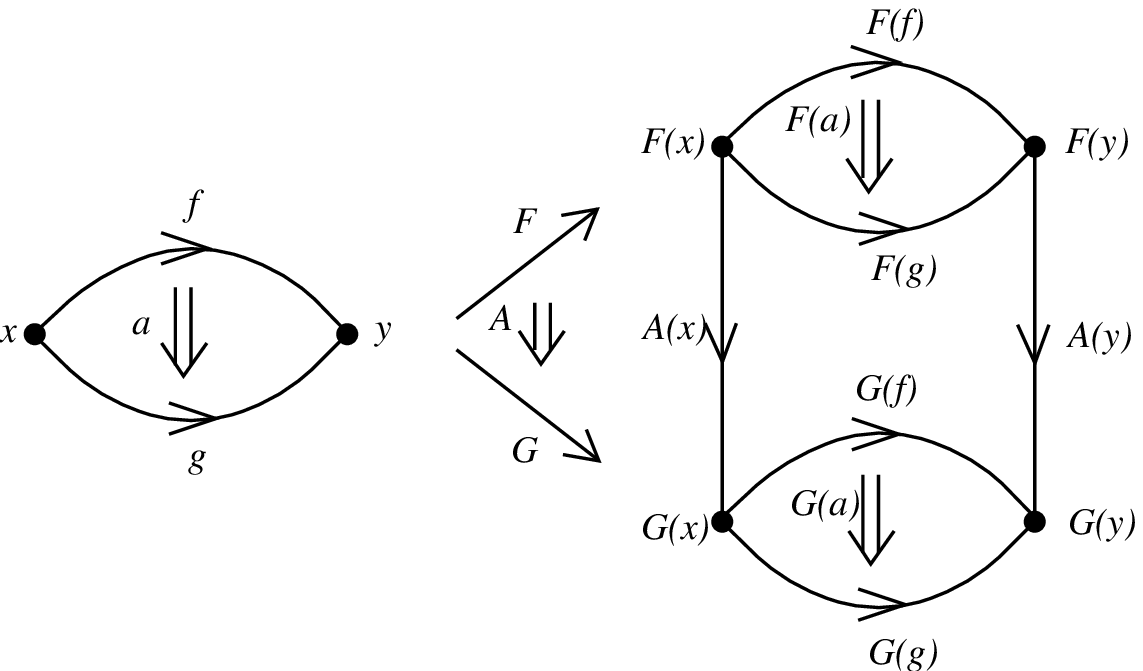}}
\medskip

\noindent Finally, the 3-morphisms are called `modifications'.  A
modification $M$ from a natural transformation $A \maps F \To G$ to a
natural transformation $B \maps F \To G$ sends each object $x \in C$ to
a 2-morphism $M(x) \maps A(x) \To B(x)$ in $D$, in a manner satisfying
some naturality conditions.  We can visualize a modification $M$ as
follows:

\medskip
\centerline{\epsfysize=2.2in\epsfbox{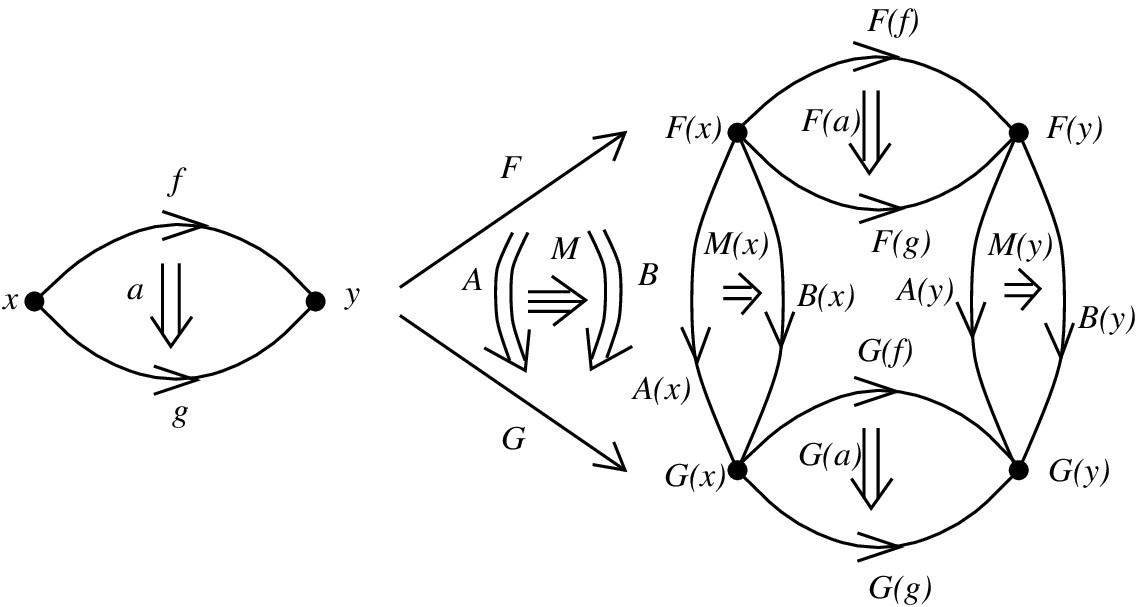}}
\medskip

\noindent Note how the $n$-dimensionality of an $n$-category leads
naturally to the $(n+1)$-dimensionality of $n$Cat.  

If instead we adopt the weak approach, we can form a tricategory Bicat
whose objects are bicategories \cite{GPS}.  The morphisms in Bicat are
called `pseudofunctors'.  A pseudofunctor $F \maps C \to D$ need not
strictly preserve all the structure.  For example, given morphisms $f
\maps x \to y$ and $g \maps y \to z$ in $C$, we do not require that
$F(fg) = F(f)F(g)$.  Instead, we require only that the two sides are
isomorphic by a specified 2-morphism, which in turn must satisfy some
coherence laws.  The 2-morphisms in Bicat are called `pseudonatural
transformations'.  In these, the squares that had to commute in the
definition of a natural transformation need only commute up to a
2-isomorphism satisfying certain coherence laws.  The 3-morphisms in
Bicat are called `modifications'.  Here there is no room for weakening,
since a modification sends each object in $C$ to a 2-morphism in $D$,
and the only sort of laws that 2-morphisms can satisfy in a bicategory
are equational laws.

One can show that Bicat is not equivalent to a strict 3-category, so we
really need the more general notion of tricategory.  Or do we?  One
might argue that we needed tricategories only because we made the
mistake of not strictifying our bicategories.  After all, every
bicategory is equivalent to a strict 2-category.  Perhaps if we replaced
every bicategory with an equivalent strict 2-category, we could work in
the strict 3-category 2Cat and never need to think about Bicat.

Alas, while superficially plausible, this line of argument is naive.  We
have said that every bicategory $C$ is equivalent to a strict 2-category $C'$.
But what does `equivalent' mean here, precisely?  It means
`equivalent, as an object of Bicat'.  In other words, for every
bicategory $C$ there is a strict 2-category $C'$ and a pseudofunctor $F
\maps C \to C'$ that is invertible up to a pseudonatural transformation
that is invertible up to a modification that is invertible!  In practice,
therefore, the business of strictifying bicategories requires a solid
understanding of Bicat as the full-fledged tricategory it is.

There is much more to say about this subject, but the basic point is
that, like it or not, sets are weak $0$-categories, and a deep
understanding of weak $n$-categories requires an understanding of weak
$(n+1)$-categories.  For this reason mathematics has been forced, over
the last century, to climb the ladder of weak $n$-categories.  To see
how computer science is repeating this climb, try for example the paper
by Power \cite{Power3} entitled `Why tricategories?'.

The actual history of this climb is quite interesting, but the details
are quite complicated, so we content ourselves here with a
thumbnail sketch.  While the formalization of the notion of set was a
slow process, the now-standard Zermelo-Fraenkel axioms reached
their final form in 1922.  Categories were defined by Eilenberg and
MacLane later in their 1945 paper \cite{EM}.  Strict 2-categories were
developed by Ehresmann \cite{Ehresmann} by 1962, and reinvented by
Eilenberg and Kelly \cite{EK} in a paper appearing the proceedings of a
conference held in 1965.  Gray \cite{Gray} discussed Cat as a strict
2-category in the same conference proceedings, and B\'enabou's
\cite{Benabou} bicategories appeared in 1967.  Gordon, Power and
Street's definition of tricategories \cite{GPS} was published in 1995,
and about this time Trimble formulated the definition of tetracategories
\cite{Trimble}.

Subsequent work has concentrated on radically accelerating this process
by defining weak $n$-categories for all $n$ simultaneously.  Actually,
Street \cite{Street} proposed a simplicial definition of weak
$n$-categories for all $n$ in 1987, but this appears not to have been
seriously studied, perhaps in part because it came too early!  Starting
in 1995, Dolan and the author gave a definition of weak $n$-categories
using `opetopic sets' \cite{BD1,BD2}, Tamsamani gave a definition using
`multisimplicial sets' \cite{Tamsamani,Tamsamani2}, and Batanin gave
definition of globular weak $\omega$-categories
\cite{Batanin,Batanin2,Street3}.  Dolan and the author have constructed
the weak $(n+1)$-category of their $n$-categories, and Simpson
\cite{Simpson} has constructed the weak $(n+1)$-category of Tamsamani's
$n$-categories.  Now the focus is turning towards working with these
different definitions and seeing whether they are equivalent.  We return
to this last issue in the Conclusions.

\subsection{Homotopy Theory} \label{homotopy}

A less inbred application of $n$-category theory is to the branch of
algebraic topology known as homotopy theory.  In fact, many
of our basic insights into $n$-categories come from this subject.
The reason is not far to seek.  Topology concerns the category Top whose
objects are topological spaces and whose morphisms are continuous
maps.  Unfortunately, there is no useful classification of topological
spaces up to isomorphism --- an isomorphism in Top being called a
`homeomorphism'.  When topologists realized this, they
retreated to the goal of classifying spaces up to various
coarser equivalence relations.  Homotopy theory is all about
properties that are preserved by continuous deformations.  More
precisely, given spaces $X,Y \in {\rm Top}$ and maps $F,G \maps X \to
Y$, one defines a `homotopy' from $F$ to $G$ to be a map 
$H \maps [0,1] \times x \to y$ with 
\[    H(0,\cdot) = F,  \qquad H(1,\cdot) = G. \]
Homotopy theory studies properties of maps that are preserved by
homotopies.  Thus two spaces $X$ and $Y$ are `the same' for the purposes
of homotopy theory, or more precisely `homotopy equivalent', if there
are maps $F \maps X \to Y$, $G \maps Y \to X$ which are inverses {\it up
to homotopy}.  

In fact, what we have done here is made Top into a 2-category whose
objects are spaces, whose morphisms are maps between spaces, and whose
2-morphisms are homotopies between maps.  This allows us to replace the
categorical concept of isomorphism between spaces by the more flexible
2-categorical concept of equivalence.  However, work on homotopy theory
soon led to the study of `higher homotopies'.  Since a homotopy is itself
a map, the concept a homotopy between homotopies makes perfect sense,
and we may iterate this indefinitely.  This amounts to treating Top as
an $n$-category for arbitrarily large $n$, or for that matter, as an
$\omega$-category.

It is worthwhile pondering how the seemingly innocuous category Top
became an $\omega$-category.  The key trick was to use the unit interval
$[0,1]$ to define higher morphisms.  The reason this trick works is that
the unit interval resembles an {\it arrow} going from $0$ to $1$.  One
could say that the abstract arrow we use in category theory is a kind of
metaphor for the unit interval --- or conversely, that the unit interval
we use in topology is a kind of metaphor for the process of going from
`here' to `there'.  However, unlike the most general sort of abstract
arrow, the unit interval has a special feature: we can go from $1$ to
$0$ as easily as we can go from $0$ to $1$.

Taking advantage of this insight, Grothendieck \cite{Gro} proposed
thinking of homotopy theory as a branch of $n$-category theory, as
follows.  We should be able to associate to any space $X$ a weak
$\omega$-category $\Pi(X)$ whose objects are points $x \in X$, whose
morphisms are paths (maps $F \maps [0,1] \to X$) going from one
point to another, whose 2-morphisms are certain paths of paths, and so
on.  Due to the special feature of the unit interval, every
$j$-morphism in this $\omega$-category should be an equivalence.
We call this special sort of $\omega$-category an `$\omega$-groupoid',
since a category with all morphisms invertible is called a groupoid.

Grothendieck also argued that conversely, we should be able to obtain a
topological space $N(G)$ from any weak $\omega$-groupoid $G$,
essentially by taking seriously the picture we can draw with points for
objects of $G$, intervals for morphisms of $G$, and so on.  By this
means we should be able to obtain weak $\omega$-functors $\Pi \maps {\rm
Top} \to \omega{\rm Gpd}$ and $N \maps \omega{\rm Gpd} \to {\rm Top}$.
Using these, we should be able to show that that the weak
$\omega$-categories Top and $\omega{\rm Gpd}$ are equivalent, as objects
of $\omega{\rm Cat}$.  In short, homotopy theory is another word for the
study of $\omega$-groupoids!

There many ways to try to realize this program, a number of which have
already obtained results.  It is well-known that all of homotopy theory
can be done purely combinatorially using `Kan complexes' \cite{May},
which may be regarded as simplicial weak $\omega$-categories.  Brown,
Higgins, Loday, and collaborators have developed a variety of approaches
using cubes \cite{Brown}.  Kapranov and Voevodsky \cite{KVinfinity} have
shown that homotopy theory is in principle equivalent to the study of
their `$\infty$-groupoids'.  Tamsamani has also shown that his approach to
weak $n$-categories reduces to homotopy theory in the $n$-groupoid case
\cite{Tamsamani2}.

Many homotopy theorists might doubt the importance of seeing homotopy
theory as a branch of $n$-category theory.  In a sense, they already
implicitly know many of the lessons $n$-category theory has to offer:
the idea of replacing equations by equivalences, the importance of
`homotopies between homotopies', and the crucial importance of coherence
laws.  Eventually $n$-category should be able to help homotopy theory in
its treatment of morphisms that are not equivalences.  In the short
term, however, the question is not what $n$-categories can do for
homotopy theory, but what homotopy theory can do for $n$-categories.

In fact, many ideas in $n$-category theory have already had their origin
in homotopy theory.  A good example is Stasheff's work on the
associahedron \cite{Stasheff}.  Recall that a `monoid' is a set equipped
with an associative product and multiplicative unit, while a
`topological monoid' is a monoid equipped with a topology for which the
product is continuous.  Stasheff wanted to uncover the
homotopy-invariant structure contained in a topological monoid.  Suppose
$X$ is a topological monoid and $Y$ is a space equipped with a homotopy
equivalence to $X$.  What sort of structure does $Y$ inherit from $X$?
Clearly we can use the homotopy equivalence to transport the product and
unit from $X$ to $Y$, obtaining a product and unit on $Y$ satisfying the
laws of a monoid {\it up to homotopy}.  For example, the two maps
\[       F,G \maps Y \times Y \times Y \to Y \]
given by 
\ban             F(y_1,y_2,y_3) &=& (y_1y_2)y_3   \\
                 G(y_1,y_2,y_3) &=& y_1(y_2y_3) \ean
need not be equal, but there is a homotopy between them, the
`associator'.  Stasheff showed that this associator satisfies the
pentagon identity up to homotopy, and that this homotopy satisfies a
coherence law of its own, again up to homotopy, and so on ad infinitum.
By working out these coherence laws in detail, he discovered the
associahedron.  Later the associahedron turned out to be relevant to
weak $\omega$-categories in general.  Part of the reason is that we can
think of the space $Y$ above as a special sort of $\omega$-category.  A
monoid can be thought of as a category with one object, by viewing the
monoid elements as morphisms from this object to itself.  Similarly, we
can view $Y$ as a weak $\omega$-category with one object, points of $Y$ as
morphisms from this object to itself, paths between these as
2-morphisms, and so on.

\subsection{Topological Quantum Field Theory}

In physics, interest in $n$-categories was sparked by developments in
relating topology and quantum field theory \cite{Kohno}.  One can
roughly date the beginning of this story to 1985, when Jones came across
a wholly unexpected invariant of knots while studying some operator
algebras invented by von Neumann in his work on the mathematical
foundations of quantum theory.  Soon this `Jones polynomial' was
generalized to a family of knot invariants.  It was then realized that
these generalizations could be systematically derived from algebraic
structures known as `quantum groups', first invented by Drinfel'd and
collaborators in their work on exactly soluble 2-dimensional field
theories \cite{CP,Majid}.  This relationship involved 2-dimensional
pictures of knots.  The story became even more exciting when Witten came
up with a manifestly 3-dimensional approach to the new knot invariants,
deriving them from a quantum field theory in 3-dimensional spacetime now
known as Chern-Simons-Witten theory.  This approach also gave invariants
of 3-dimensional manifolds.

These developments exposed a deep but mysterious unity in what at first
might seem like disparate branches of algebra, topology, and quantum
physics.  Interestingly, it appears that the roots of this unity lie in
certain aspects of $n$-category theory.  In fact, this is the main
reason for the author's interest in $n$-categories: it appears that a
good theory of weak $n$-categories is necessary to serve as a framework
for the mathematics that will be able to reduce the currently rather
elaborate subject of `topological quantum field theory' to its simple
essence.  Having explained this at length elsewhere \cite{BD}, we limit
our remarks here to a few key points.  

Quantum physics relies crucially on the theory of Hilbert spaces.  For
simplicity, we limit our attention here to the finite-dimensional case,
defining a `Hilbert space' to be a finite-dimensional complex vector
space $H$ equipped with an `inner product'
\[            \langle \cdot, \cdot \rangle \maps H \times H \to \C \]
which is linear in the second argument, conjugate-linear in the first,
and satisfies $\langle \psi, \phi \rangle = \overline {\langle \phi,
\psi \rangle}$ for all $\psi,\phi \in H$ and $\langle \psi,\psi \rangle
> 0$ for all nonzero $\psi \in H$.  The inner product allows us to
define the norm of a vector $\psi \in H$ by
\[       \|\psi \| = \langle \psi , \psi \rangle^{1/2} ,\]
but its main role in physics to compute amplitudes.  States of a quantum
system are described by vectors with norm 1.  If one places a quantum
system in the state $\psi$, and then does an experiment to see if it is
in some state $\phi$, the probability that the answer is `yes' equals
\[              |\langle \phi,\psi \rangle|^2  .\]
This is automatically a real number between 0 and 1.  However, when one
delves deeper into the theory, it appears that even more fundamental
than the probability is the `amplitude'  
\[          \langle \phi, \psi \rangle , \]
which is of course a complex number.  

The role of the inner product in quantum physics has always been a
source of puzzles to those with an interest in the philosophical
foundations of the subject.  Complex amplitudes lack the intuitive
immediacy of probabilities.  From the category-theoretic point of view,
part of the problem is to understand the category of Hilbert spaces.
The objects are Hilbert spaces, but what are the morphisms?  Typically
morphisms are required to preserve all the structure in sight.  This
suggests taking the morphisms to be linear operators preserving the
inner product.  However, other linear operators are also important.
Particularly in topological quantum field theory, there are good reasons
to take {\it all} linear operators as morphisms.  However, if we define
a category Hilb this way, then Hilb is equivalent to the category Vect
of complex vector spaces.  This then raises the question: how does Hilb
really differ from Vect, if as categories they are equivalent?

Luckily, quantum theory suggests an answer to this question.  Given any 
linear operator $F \maps H \to H'$ between Hilbert spaces, we may define the
`adjoint' $F^\ast \maps H' \to H$ to be the unique linear operator with
\[    \langle F\phi,\psi \rangle = \langle \phi,F^\ast \psi\rangle \]
for all $\phi \in H$, $\psi \in H'$.  This sort of adjoint is basic to
quantum theory.  From the category-theoretic point of view, the role of
the adjoint is to make Hilb into a `$\ast$-category': a category $C$
equipped with a contravariant functor $\ast \maps C \to C$ fixing
objects and satisfying $\ast^2 = 1_C$.  While Hilb and Vect are
equivalent as categories, only Hilb is a $\ast$-category.

This is particularly important in topological quantum field theory.  The
mysterious relationships between topology, algebra and physics exploited
by this subject amount in large part to the existence of interesting
functors from various topologically defined categories to the category
Hilb.  These topologically defined categories are always
$\ast$-categories, and the really interesting functors from them to Hilb
are always `$\ast$-functors', functors preserving the $\ast$-structure.
Physically, the $\ast$ operation corresponds to {\it reversing the
direction of time}.   For example, there is a $\ast$-category whose
objects are collections of points and whose morphisms are `tangles':

\medskip
\centerline{\epsfysize=1.5in\epsfbox{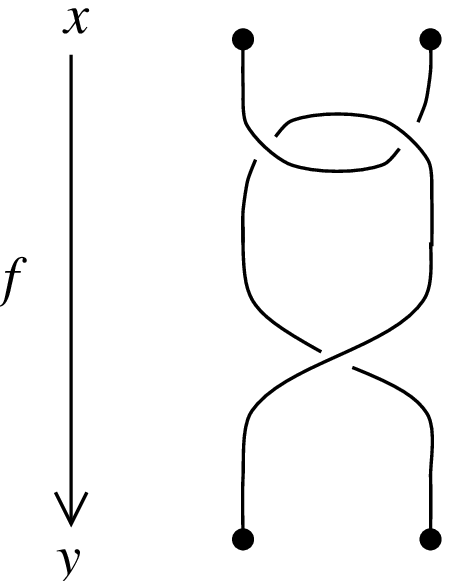}}
\medskip

\noindent We can think of this morphism $f \maps x \to y$ as
representing the trajectories of a collection of particles and
antiparticles, where particles and antiparticles can be created or
annihilated in pairs.  Reversing the direction of time, we obtain the
`dual' morphism $f^\ast \maps y \to x$:

\medskip
\centerline{\epsfysize=1.5in\epsfbox{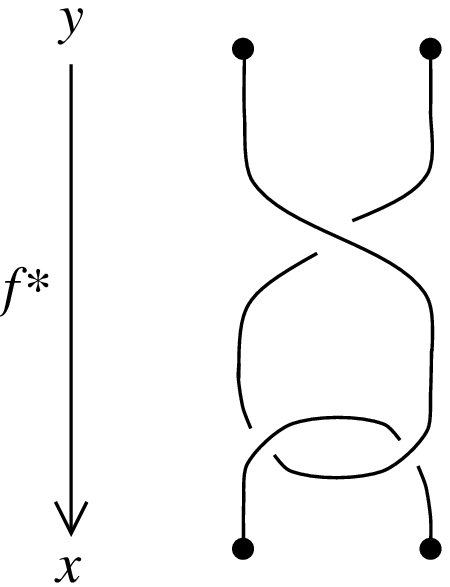}}
\medskip

\noindent This morphism is not the inverse of $f$, since the composite
$ff^\ast$ is a nontrivial tangle:

\medskip\medskip
\centerline{\epsfysize=3.0in\epsfbox{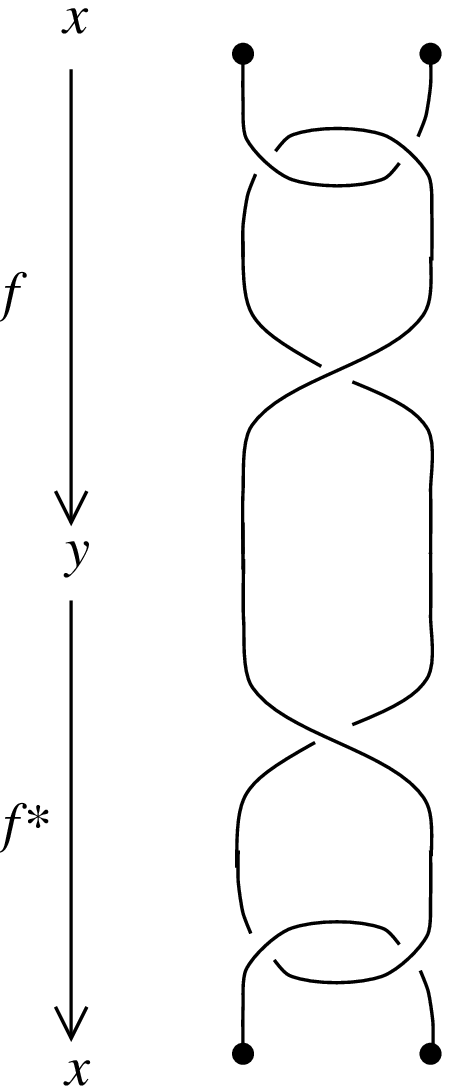}}
\medskip

\noindent Indeed, any groupoid becomes a $\ast$-category if we set
$f^\ast = f^{-1}$ for every morphism $f$, but the most interesting
$\ast$-categories in topological quantum field theory are not groupoids.

The above example involves 1-dimensional curves in 3-dimensional
spacetime.  More generally, topological quantum field theory studies
$n$-dimensional manifolds embedded in $(n+k)$-dimensional spacetime,
which in the $k \to \infty$ limit appear as `abstract' $n$-dimensional
manifolds.   It appears that these are best described using certain
`$n$-categories with duals', meaning $n$-categories in which 
every $j$-morphism $f$ has a dual $f^\ast$.   Unfortunately, so far
the details have only been worked out in certain low-dimensional cases
\cite{B,BL}.  The main problem is that the notion of `$n$-category with
duals' is only beginning to be understood.  

One class of $n$-categories with duals should be the $n$-groupoids; this
would explain many relationships between topological quantum field
theory and homotopy theory \cite{Porter}.  However, the novel aspects of
topological quantum field theory should arise from $n$-categories with
duals that are not $n$-groupoids.  Indeed, this explains why the Jones
polynomial and other new knot invariants were not discovered earlier
using traditional techniques of algebraic topology.

The idea that duals are subtler and thus more interesting than inverses
is already familiar from category theory.  Given a functor $F
\maps C \to D$, the correct sort of weakened `inverse' to $F$ is a 
functor $G \maps D \to C$ such that $FG$ and $GF$ are naturally
isomorphic to the identity; if such a functor $G$ exists then $F$ is an
equivalence.  However, even if no such `inverse' exists, the functor $F$
may have a kind of `dual', namely an adjoint functor!  A right adjoint
$F^\ast \maps C \to D$, for example, would satisfy:
\[    \hom(Fx,y) \iso \hom(x,F^\ast y) \]
for all $x \in C$, $y \in D$.  Note that this is very similar to the
definition of the adjoint of a linear map between Hilbert spaces, with
`hom' playing the role of the inner product.  

The analogy between adjoint functors and adjoint linear operators relies
upon a deeper analogy: just as in quantum theory the inner product
$\langle \phi,\psi\rangle$ represents the {\it amplitude} to pass from
$\phi$ to $\psi$, in category theory $\hom(x,y)$ represents the {\it set
of ways} to go from $x$ to $y$.  A precise working out of this analogy
can be found in the author's paper \cite{B} on `2-Hilbert spaces'.
These are to Hilbert spaces as categories are to sets.  The analogues of 
adjoint linear operators between Hilbert spaces are certain adjoint
functors between 2-Hilbert spaces.  Just as the primordial example of a
category is Set, the primordial example of a 2-Hilbert space is Hilb.
Also, just as the 2-category Cat is a 3-category, it appears that the
2-category 2Hilb is an example of a `3-Hilbert space' --- a concept
which has not yet been given a proper definition.  

More generally, it appears that $n$Hilb is an $n$-category with duals,
and that `$n$-Hilbert spaces' are needed for the proper treatment of
$n$-dimensional topological quantum field theories \cite{BD,Freed}.
Thus, just as mathematics has been forced to ascend the ladder of
$n$-categories, so may be physics!

\section{A Definition of Weak $n$-Category}

As discussed in Section 2, any definition of $n$-categories involves a
choice of the basic shapes of $j$-morphisms and a choice of allowed ways
to glue them together.  Any definition of weak $n$-categories also
requires a careful treatment of coherence laws.  In what follows we
present a definition of weak $n$-categories in which all these issues
are handled in a tightly linked way.  In this definition, the basic
shapes of $j$-morphisms are the $j$-dimensional `opetopes'.  The allowed
ways of gluing together the $j$-dimensional opetopes correspond
precisely to the $(j+1)$-dimensional opetopes.  Moreover, the coherence
laws satisfied by composition correspond to still higher-dimensional
opetopes!

Before going into the details, let us give a rough sketch of this works.
First consider some low-dimensional opetopes.  The only 0-dimensional
opetope is the point:

\medskip
\centerline{\epsfysize=0.08in\epsfbox{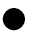}}
\medskip

\noindent There is no way to glue together 0-dimensional opetopes.  The
only 1-dimensional opetope is the interval, or more precisely the arrow:

\medskip
\centerline{\epsfysize=0.12in\epsfbox{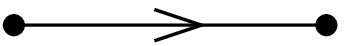}}
\medskip

\noindent The allowed ways of gluing together 1-dimensional opetopes are
given by the 2-dimensional opetopes.  The first few 2-dimensional
opetopes are as follows:

\medskip
\centerline{\epsfysize=0.8in\epsfbox{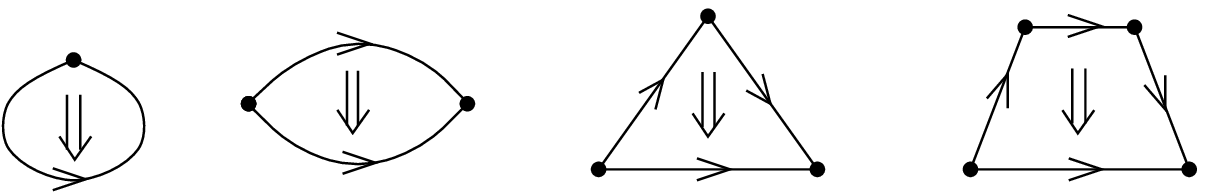}}
\medskip
 
\noindent For any $k \ge 0$, there is a 2-dimensional opetope with $k$
`infaces' and one `outface'.  (We are glossing over some subleties here;
for reasons noted later, there are really $k!$ such opetopes.)

The allowed ways of gluing together 2-dimensional opetopes are
given by the 3-dimensional opetopes.  There are many of these; a simple
example is as follows:

\medskip
\centerline{\epsfysize=1.2in\epsfbox{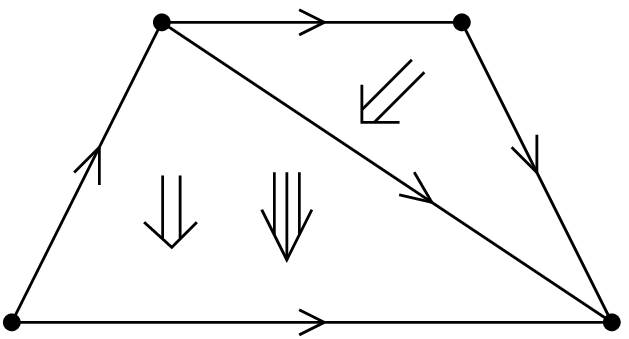}}
\medskip

\noindent This may be a bit hard to visualize, but it depicts a
3-dimensional shape whose front consists of two 3-sided `infaces', and
whose back consists of a single 4-sided `outface'.  We have drawn double
arrows on the infaces but not on the outface.  Note that while this shape is
topologically a ball, it cannot be realized as a polyhedron with
planar faces.  This is typical of opetopes.   

In general, a $(n+1)$-dimensional opetope has any number of infaces and
exactly one outface: the infaces are $n$-dimensional opetopes glued
together in a tree-like pattern, while the outface is a single
$n$-dimensional opetope.  For example, the 3-dimensional opetope above
corresponds to the following tree:

\medskip
\centerline{\epsfysize=1.2in\epsfbox{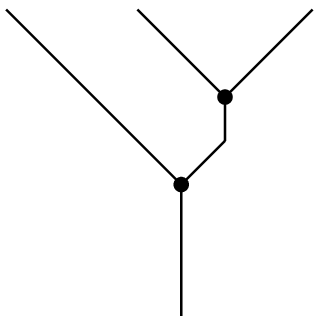}}
\medskip

\noindent The two triangular infaces of the opetope correspond to the
two nodes in this tree.  This is a rather special tree; in general, we
allow nonplanar trees with any number of nodes and any number of edges
coming into each node.

In our approach, a weak $n$-category is a special sort of `opetopic set'.
Basically, an opetopic set is a set of `cells' shaped like opetopes,
such that any face of any cell is again a cell.  In a weak $n$-category,
the $j$-dimensional cells play the role of $j$-morphisms.  An opetopic
set $C$ is an $n$-category if it satisfies the following two properties:

1) {\it ``Any niche has a universal occupant.''}  A `niche' is
a configuration where all the infaces of an opetope have been filled in
by cells of $C$, but not the outface or the opetope itself:

\medskip
\centerline{\epsfysize=1.5in\epsfbox{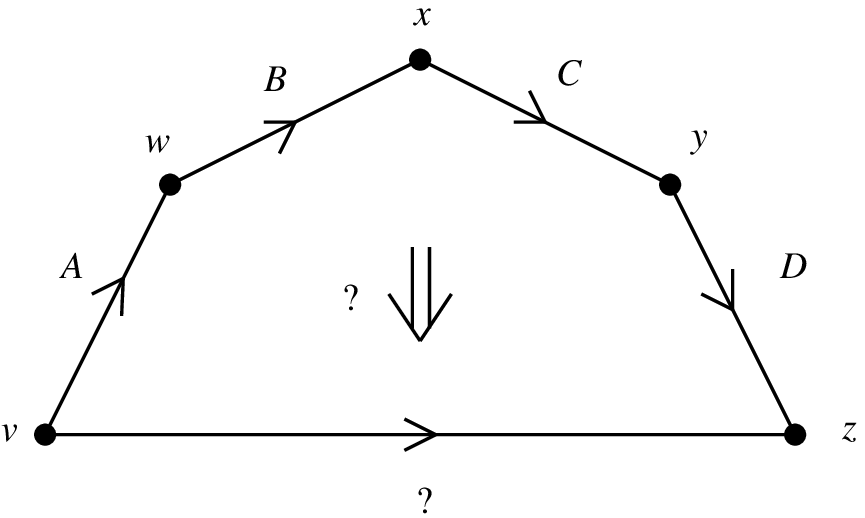}}
\medskip

\noindent An `occupant' of the niche is a way of extending this configuration
by filling in the opetope (and thus its ouface) with a cell:

\medskip
\centerline{\epsfysize=1.5in\epsfbox{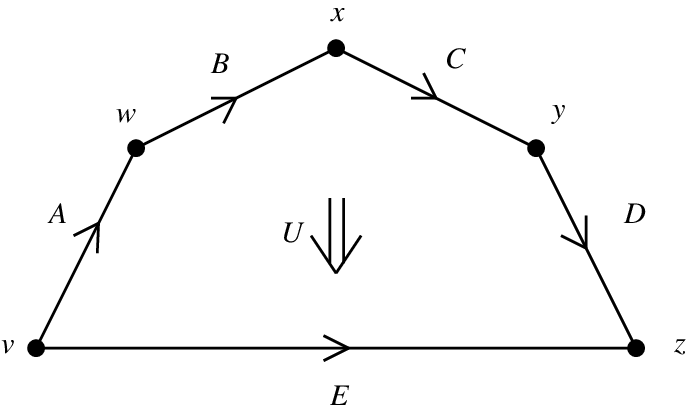}}
\medskip

\noindent The `universality' of an occupant means roughly that every other
occupant factors through the given one {\it up to equivalence}.  To make
this precise we need to define universality in a rather subtle recursive
way.  We may think of a universal occupant of a niche as `a process of
composing' the infaces, and its outface as `a composite' of the infaces.

2) {\it ``Composites of universal cells are universal.''}  Suppose
that $U,V,$ and $W$ below are universal cells:

\medskip
\centerline{\epsfysize=2.0in\epsfbox{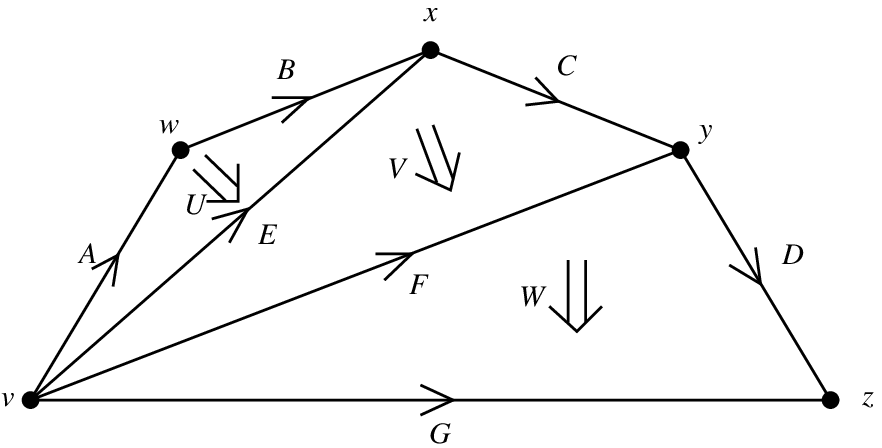}}
\medskip

\noindent Then we can compose them, and we are guaranteed that their
composite is again universal, and thus that the outface $G$ is a composite
of the cells $A,B,C,D$.  Note that a process of composing $U,V,W$
is described by a universal occupant of a niche of one higher
dimension.  

Note that in this approach to weak $n$-categories, composition of cells
is not an operation in the traditional sense: the composite is defined
by a universal property, and is thus unique only {\it up to
equivalence}.  Only at the top level, for the $n$-cells in an
$n$-category, is the composite truly unique.  This may seem odd at first
glance, but in fact it closely reflects actual mathematical practice.
For example, it is unnatural to think of the Cartesian product as an
operation on Set.  We can do it, but there are as many ways to do so as
there are ways to define ordered pairs in set theory; there is certainly
nothing sacred about the von Neumann definition.  If we arbitrarily
choose a way, we can think of Set as a weak monoidal category, i.e., a
bicategory with one object.  However, we can avoid this arbitrariness if
we define the Cartesian product of sets by a universal property, using
the category-theoretic concept of `product'.  Then the product of sets
is only defined up to a natural isomorphism, and Set becomes our sort of
weak 2-category with one object.  In this approach, all the necessary 
coherence laws {\it follow automatically} from the universal property
defining the product.  

In the following sections we first review the theory of operads,
and then explain how this theory can be used to define the opetopes.
After a brief discussion of some notions related to opetopic sets, we
give the definition of `universal occupant' of a niche, and then the
definition of weak $n$-category.  At various points we skim over
technical details; these can all be found in our paper \cite{BD2}.

\subsection{Operads}

To describe the opetopes we need to specify exactly which tree-like
patterns we can use to glue together the opetopes in one dimension to
form an opetope of the next higher dimension.  For this we use the
theory of `operads'.  An operad is a gadget consisting of abstract
$k$-ary `operations' and various ways of composing them, and the
$n$-dimensional opetopes will be the operations of a certain operad.
This is another example of how $n$-category theory is indebted to
homotopy theory, since operads were first developed for the purposes of
homotopy theory \cite{Adams,BV,May}.

In what follows we work with typed operads having a set $S$ of types,
or `$S$-operads' for short.  The basic idea of an $S$-operad $O$ is that
given types $x_1,\dots,x_k,x' \in S$, there is a set
$O(x_1,\dots,x_k;x')$ of abstract $k$-ary `operations' with inputs of
type $x_1,\dots,x_k$ and output of type $x'$.  We can visualize such an
operation as a tree with only one node:

\begin{center}
\setlength{\unitlength}{0.0007500in}%
\begingroup\makeatletter\ifx\SetFigFont\undefined
\def\x#1#2#3#4#5#6#7\relax{\def\x{#1#2#3#4#5#6}}%
\expandafter\x\fmtname xxxxxx\relax \def\y{splain}%
\ifx\x\y   
\gdef\SetFigFont#1#2#3{%
  \ifnum #1<17\tiny\else \ifnum #1<20\small\else
  \ifnum #1<24\normalsize\else \ifnum #1<29\large\else
  \ifnum #1<34\Large\else \ifnum #1<41\LARGE\else
     \huge\fi\fi\fi\fi\fi\fi
  \csname #3\endcsname}%
\else
\gdef\SetFigFont#1#2#3{\begingroup
  \count@#1\relax \ifnum 25<\count@\count@25\fi
  \def\x{\endgroup\@setsize\SetFigFont{#2pt}}%
  \expandafter\x
    \csname \romannumeral\the\count@ pt\expandafter\endcsname
    \csname @\romannumeral\the\count@ pt\endcsname
  \csname #3\endcsname}%
\fi
\fi\endgroup
\begin{picture}(1897,2565)(3226,-3304)
\thicklines
\put(4201,-1861){\circle*{150}}
\put(4201,-3061){\line( 0, 1){1200}}
\put(4201,-1861){\line(-1, 3){300}}
\put(4201,-1861){\line( 1, 3){300}}
\put(4201,-1861){\line( 1, 1){900}}
\put(4201,-1861){\line(-1, 1){900}}
\put(3226,-811){\makebox(0,0)[lb]{\raisebox{0pt}[0pt][0pt]{$x_1$}}}
\put(4126,-3286){\makebox(0,0)[lb]{\raisebox{0pt}[0pt][0pt]{$x'$}}}
\put(3901,-811){\makebox(0,0)[lb]{\raisebox{0pt}[0pt][0pt]{$x_2$}}}
\put(4501,-811){\makebox(0,0)[lb]{\raisebox{0pt}[0pt][0pt]{$x_3$}}}
\put(5101,-811){\makebox(0,0)[lb]{\raisebox{0pt}[0pt][0pt]{$x_4$}}}
\put(4426,-1966){\makebox(0,0)[lb]{\raisebox{0pt}[0pt][0pt]{$f$}}}
\end{picture}
\end{center}
\medskip
\centerline{An operation $f \in O(x_1,\dots,x_k;x')$ }
\medskip

\noindent In an operad, we can get new operations
from old ones by composing them, which we can visualize 
in terms of trees as follows:

\begin{center}
\setlength{\unitlength}{0.000500in}%
\begingroup\makeatletter\ifx\SetFigFont\undefined
\def\x#1#2#3#4#5#6#7\relax{\def\x{#1#2#3#4#5#6}}%
\expandafter\x\fmtname xxxxxx\relax \def\y{splain}%
\ifx\x\y   
\gdef\SetFigFont#1#2#3{%
  \ifnum #1<17\tiny\else \ifnum #1<20\small\else
  \ifnum #1<24\normalsize\else \ifnum #1<29\large\else
  \ifnum #1<34\Large\else \ifnum #1<41\LARGE\else
     \huge\fi\fi\fi\fi\fi\fi
  \csname #3\endcsname}%
\else
\gdef\SetFigFont#1#2#3{\begingroup
  \count@#1\relax \ifnum 25<\count@\count@25\fi
  \def\x{\endgroup\@setsize\SetFigFont{#2pt}}%
  \expandafter\x
    \csname \romannumeral\the\count@ pt\expandafter\endcsname
    \csname @\romannumeral\the\count@ pt\endcsname
  \csname #3\endcsname}%
\fi
\fi\endgroup
\begin{picture}(3405,4252)(3879,-5183)
\thicklines
\put(4201,-1861){\circle*{150}}
\put(5701,-1861){\circle*{150}}
\put(7201,-1861){\circle*{150}}
\put(5701,-3661){\circle*{150}}
\put(4201,-3061){\line( 0, 1){1200}}
\put(4201,-3061){\line( 5,-2){1500}}
\put(5701,-3661){\line( 0, 1){600}}
\put(5701,-3061){\line( 0, 1){1200}}
\put(5701,-3661){\line( 5, 2){1500}}
\put(7201,-3061){\line( 0, 1){1200}}
\put(7201,-1861){\line( 0, 1){900}}
\put(5701,-3661){\line( 0,-1){1500}}
\put(5701,-1861){\line( 1, 1){900}}
\put(5678,-1853){\line( 1, 3){300}}
\put(5723,-1853){\line(-1, 3){300}}
\put(5701,-1861){\line(-1, 1){900}}
\put(4201,-1861){\line( 1, 3){300}}
\put(4201,-1936){\line(-1, 3){322.500}}
\put(3851,-1936){\makebox(0,0)[lb]{\raisebox{0pt}[0pt][0pt]{$g_1$}}}
\put(5351,-1936){\makebox(0,0)[lb]{\raisebox{0pt}[0pt][0pt]{$g_2$}}}
\put(6851,-1936){\makebox(0,0)[lb]{\raisebox{0pt}[0pt][0pt]{$g_3$}}}
\put(5326,-3886){\makebox(0,0)[lb]{\raisebox{0pt}[0pt][0pt]{$f$}}}
\end{picture}
\end{center}
\medskip
\centerline{An operation $f \circ (g_1,\dots,g_k)$}
\medskip

\noindent 
We can also obtain new operations from old by permuting arguments, and
there is a unary `identity' operation of each type.  Finally, we demand 
a few plausible axioms: the identity operations act as identities for
composition, permuting arguments is compatible with composition, and
composition is `associative', making composites of the following sort
well-defined:

\vbox{

\begin{center}

\setlength{\unitlength}{0.000250in}%
\begingroup\makeatletter\ifx\SetFigFont\undefined
\def\x#1#2#3#4#5#6#7\relax{\def\x{#1#2#3#4#5#6}}%
\expandafter\x\fmtname xxxxxx\relax \def\y{splain}%
\ifx\x\y   
\gdef\SetFigFont#1#2#3{%
  \ifnum #1<17\tiny\else \ifnum #1<20\small\else
  \ifnum #1<24\normalsize\else \ifnum #1<29\large\else
  \ifnum #1<34\Large\else \ifnum #1<41\LARGE\else
     \huge\fi\fi\fi\fi\fi\fi
  \csname #3\endcsname}%
\else
\gdef\SetFigFont#1#2#3{\begingroup
  \count@#1\relax \ifnum 25<\count@\count@25\fi
  \def\x{\endgroup\@setsize\SetFigFont{#2pt}}%
  \expandafter\x
    \csname \romannumeral\the\count@ pt\expandafter\endcsname
    \csname @\romannumeral\the\count@ pt\endcsname
  \csname #3\endcsname}%
\fi
\fi\endgroup
\begin{picture}(4980,5444)(3218,-5183)
\thicklines
\put(4201,-1861){\circle*{150}}
\put(5701,-1861){\circle*{150}}
\put(7201,-1861){\circle*{150}}
\put(5701,-3661){\circle*{150}}
\put(3901,-361){\circle*{150}}
\put(4501,-361){\circle*{150}}
\put(5101,-361){\circle*{150}}
\put(6601,-361){\circle*{150}}
\put(7801,-361){\circle*{150}}
\put(3301,-361){\circle*{150}}
\put(7201,-361){\circle*{150}}
\put(4201,-3061){\line( 0, 1){1200}}
\put(4201,-1861){\line(-1, 3){300}}
\put(4201,-1861){\line( 1, 3){300}}
\put(4201,-1861){\line( 1, 1){900}}
\put(4201,-1861){\line(-1, 1){900}}
\put(4201,-3061){\line( 5,-2){1500}}
\put(5701,-3661){\line( 0, 1){600}}
\put(5701,-3061){\line( 0, 1){1200}}
\put(5701,-3661){\line( 5, 2){1500}}
\put(7201,-3061){\line( 0, 1){1200}}
\put(7201,-1861){\line(-2, 3){600}}
\put(7201,-1861){\line( 0, 1){900}}
\put(7201,-1861){\line( 2, 3){600}}
\put(5701,-3661){\line( 0,-1){1500}}
\put(3301,-961){\line( 0, 1){600}}
\put(3901,-961){\line( 0, 1){600}}
\put(4501,-961){\line( 0, 1){600}}
\put(5101,-961){\line( 0, 1){600}}
\put(6601,-961){\line( 0, 1){600}}
\put(7201,-1036){\line( 0, 1){675}}
\put(7801,-961){\line( 0, 1){600}}
\put(5101,-361){\line( 1, 4){150}}
\put(5101,-361){\line(-1, 4){150}}
\put(4501,-361){\line( 0, 1){600}}
\put(3901,-361){\line( 1, 2){300}}
\put(3901,-361){\line( 0, 1){600}}
\put(3901,-361){\line(-1, 2){300}}
\put(6601,-361){\line( 1, 2){300}}
\put(6601,-361){\line( 0, 1){600}}
\put(6601,-436){\line(-2, 5){274.138}}
\put(7801,-361){\line(-3, 5){363.971}}
\put(7801,-361){\line(-1, 4){150}}
\put(7801,-361){\line( 1, 4){150}}
\put(7801,-361){\line( 3, 5){363.971}}
\end{picture}
\end{center}
\medskip
}

Formally, we have:

\begin{defn}\et For any set $S$, an `$S$-operad' $O$ consists of
{\rm \begin{enumerate}  
\item {\it for any $x_1,\dots,x_k,x' \in S$, a set
$O(x_1,\dots,x_k;x')$}
\item  {\it for any $f \in O(x_1,\dots,x_k;x')$ and any
$g_1 \in O(x_{11},\dots,x_{1i_1};x_1), \dots,$ \hfill \break 
$g_k \in O(x_{k1},\dots,x_{ii_k};x_k)$, an element
\[ f \circ (g_1,\dots, g_k) \in O(x_{11},\dots,x_{1i_1}, \dots \dots,
x_{k1}, \dots , x_{ii_k};x') \] }
\item {\it for any $x \in S$, an element $1_x \in O(x;x)$}
\item {\it for any permutation $\sigma \in S_k$, a map
\ban     \sigma \maps O(x_1,\dots,x_k;x') &\to&
          O(x_{\sigma(1)}, \dots, x_{\sigma(k)};x')  \\
f &\mapsto& f\sigma \ean  }
\end{enumerate}}
\noindent such that:
\begin{alphalist}
\item  whenever both sides make sense, 
\[      f \cdot (g_1 \cdot (h_{11}, \dots, h_{1i_1}), \dots, 
g_k \cdot (h_{k1}, \dots, h_{ki_k})) = \]
\[ (f \cdot (g_1, \dots g_k)) \; \cdot \; (h_{11}, \dots, h_{1i_1},
\dots\dots , h_{k1}, \dots, h_{ki_k}) \]
\item  for any $f \in O(x_1,\dots,x_k;x')$,
\[ f =  1_{x'} \cdot f = f \cdot (1_{x_1},\dots,1_{x_k}) \]
\item  for any $f \in O(x_1,\dots,x_k;x')$ and 
$\sigma,\sigma' \in S_k$, 
\[      f(\sigma \sigma') = (f \sigma)\sigma' \]
\item  for any $f \in O(x_1,\dots,x_k;x')$, $\sigma \in S_k$, and 
$g_1 \in O(x_{11},\dots,x_{1i_1};x_1),\dots,$ \hfill \break
$g_k \in O(x_{k1},\dots,x_{ki_k};x_k)$, 
\[     (f\sigma) \cdot (g_{\sigma(1)}, \dots, g_{\sigma(k)}) = (f \cdot
(g_1, \dots, g_k))\, \rho(\sigma), \]
where $\rho \maps S_k \to S_{i_1 + \cdots + i_k}$ is the obvious
homomorphism.
\item  for any $f \in O(x_1,\dots,x_k;x')$, 
$g_1 \in O(x_{11},\dots,x_{1i_1};x_1),\dots,$ \hfill \break
$g_k \in O(x_{k1},\dots,x_{ki_k};x_k)$, and 
$\sigma_1 \in S_{i_1}, \dots, \sigma_k \in S_{i_k}$,
\[     (f\cdot (g_1\sigma_1 , \dots, g_k\sigma_k))  = (f \cdot
(g_1, \dots, g_k))\, \rho'(\sigma_1,\dots,\sigma_k), \]
where $\rho' \maps S_{i_1} \times \cdots \times S_{i_k}
\to S_{i_1 + \cdots + i_k}$ is the obvious homomorphism.
\end{alphalist}
\end{defn}

Operads have `algebras', in which their abstract operations are 
represented as actual functions:

\begin{defn}  \label{c.op.algebras}\et  For any $S$-operad
$O$, an `$O$-algebra' $A$ consists of:
{\rm \begin{enumerate}
\item {\it for any $x \in S$, a set $A(x)$.}
\item {\it for any $f \in O(x_1,\dots,x_k;x')$, a function
\[   \alpha(f) \maps A(x_1) \times \cdots \times A(x_k) \to A(x') \]  }
\end{enumerate}}
\noindent such that:
\begin{alphalist}
\item whenever both sides make sense,
\[      \alpha(f \cdot (g_1, \dots, g_k)) = \alpha(f) (\alpha(g_1)
\times \cdots \times \alpha(g_k)) \]
\item  for any $x \in C$, $\alpha(1_x)$ acts as the
identity on $A(x)$
\item for any $f \in O(x_1,\dots,x_k,x')$ and $\sigma \in S_k$,
\[         \alpha(f\sigma) = \alpha(f)\sigma, \]
where $\sigma \in S_k$ acts on the function $\alpha(f)$ 
on the right by permuting its arguments. 
\end{alphalist}
\end{defn}

We can think of an operad as a simple sort of theory, and its
algebras as models of this theory.   Thus we can study operads
either `syntactically' or `semantically'.  To describe an operad
syntactically, we list:
\begin{enumerate}
\item the set $S$ of {\it types},
\item the sets $O(x_1,\dots,x_k;x')$ of {\it operations},
\item the set of all {\it reduction laws} saying that some composite of
operations (possibly with arguments permuted) equals some other
operation.  
\end{enumerate}
This is like a presentation in terms of generators and relations, with
the reduction laws playing the role of relations.  On the other hand, to
describe an operad semantically, we describe its algebras.  

\subsection{Opetopes}

The following fact is the key to defining the opetopes.  Let $O$ be an
$S$-operad, and let $\elt(O)$ be the set of all operations of $O$.

\begin{thm}\et There is an $\elt(O)$-operad $O^+$ whose algebras are
$S$-operads over $O$, i.e., $S$-operads equipped with a homomorphism to
$O$.  \end{thm}

\noindent We call $O^+$ the `slice operad' of $O$.  One can describe 
$O^+$ syntactically as follows:
\begin{enumerate}
\item The types of $O^+$ are the operations of $O$.
\item The operations of $O^+$ are the reduction laws of $O$.
\item The reduction laws of $O^+$ are the ways of combining reduction
laws of $O$ to give other reduction laws.
\end{enumerate}
The `level-shifting' going on here as we pass from $O$ to $O^+$
captures the process by which equational laws are promoted to equivalences and 
these equivalences satisfy new coherence laws of their own.  In this
context, the new laws are just {\it the ways of combining the the old laws}.

Note that we can iterate the slice operad construction. 
Let $O^{n+}$ denote the result of applying the slice operad construction
$n$ times to the operad $O$ if $n \ge 1$, or just $O$ itself if $n = 0$.

\begin{defn}\et An $n$-dimensional `$O$-opetope' is a type of $O^{n+}$,
or equivalently, if $n \ge 1$, an operation of $O^{(n-1)+}$.  \end{defn}

\noindent In particular, we define an $n$-dimensional `opetope' to be an
$n$-dimensional $O$-opetope when $O$ is the simplest operad of all:

\begin{defn}\et The `initial untyped operad' $I$ is the $S$-operad with:
{\rm \begin{enumerate}
\item {\it only one type: $S = \{x\}$}
\item {\it only one operation, the identity operation $1 \in O(x;x)$}
\item {\it all possible reduction laws}
\end{enumerate}
}
\end{defn}

\noindent Semantically, $I$ is the operad whose algebras are just sets.

The opetopes emerge from $I$ as follows.  The 0-dimensional opetopes
are the types of $I$, but there is only one type, so there is only
one 0-dimensional opetope:

\medskip
\medskip
\centerline{\epsfysize=0.2in\epsfbox{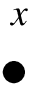}}
\medskip

\noindent The 1-dimensional opetopes are the types of $I^+$, or in other
words, the operations of $I$.  There is only one operation of $I$, the
identity operation, so there is only one 1-dimensional opetope:

\medskip
\centerline{\epsfysize=0.3in\epsfbox{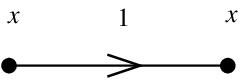}}
\medskip

\noindent The 2-dimensional opetopes are the types of $I^{++}$, or in
other words, the operations of $I^+$, which are the reduction laws of
$I$.  These reduction laws all state that the identity operation
composed with itself $k$ times equals itself.  This leads to
2-dimensional opetopes with $k$ infaces and one outface, as follows:

\medskip
\centerline{\epsfysize=1.0in\epsfbox{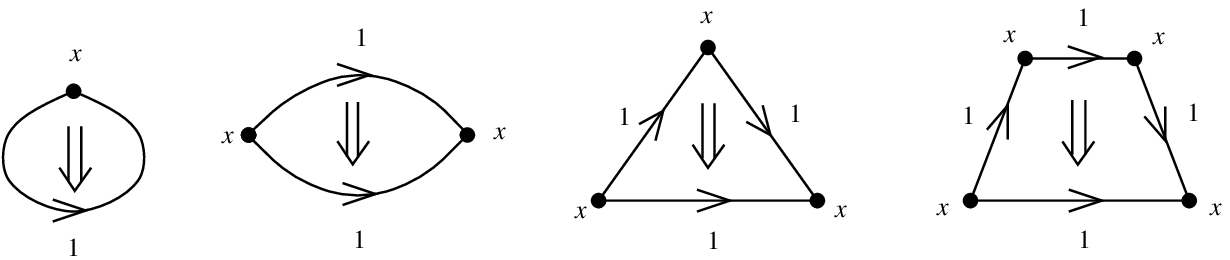}}
\medskip

\noindent Actually there are $k!$ different 2-dimensional opetopes with
$k$ infaces, since the permutation group $S_k$ acts freely on the set of
$k$-ary operations of $I^+$.  We could keep track of these by labelling
the infaces with some permutation of $k$ distinct symbols.

The 3-dimensional opetopes are the types of $I^{+++}$, or in other
words, the operations of $I^{++}$, which are the reduction laws of
$I^+$.  These state that some composite of 2-dimensional opetopes
equals some other 2-dimensional opetope.  This leads to 3-dimensional 
opetopes like the following:

\medskip
\centerline{\epsfysize=1.4in\epsfbox{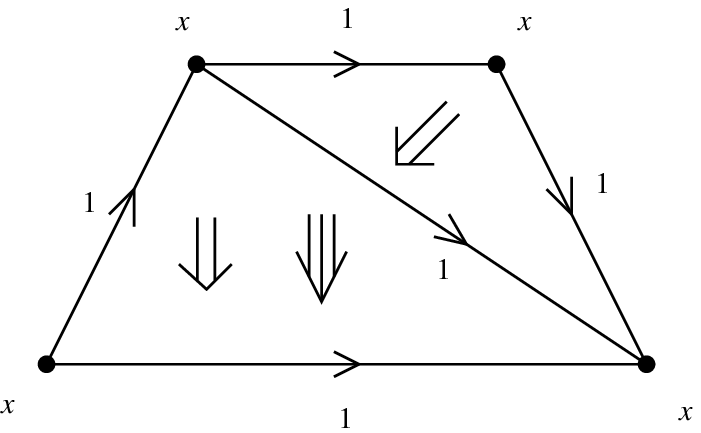}}
\medskip

In general, the $(n+1)$-dimensional opetopes describe all possible ways
of composing $n$-dimensional opetopes.  Since the $(n+1)$-dimensional
opetopes are the operations of the operad $I^{n+}$, all allowed ways of
composing them can be described by trees.  One can use this to describe
the opetopes of all dimensions using `metatree notation'.  In this
notation, an $n$-dimensional opetope is represented as a list of $n$
labelled trees.  This notation nicely handles the nuances of how the
permutation group $S_k$ acts on the set of opetopes with $k$ infaces.

\subsection{Opetopic sets}

A weak $n$-category will be an `opetopic set' with certain properties.
An opetopic set consists of collections of `cells' of different shapes,
one collection for each opetope.  The face of any cell is again a cell,
and one can keep track of this using `face maps' going between the
collections of cells.  These maps also satisfy certain relations.  We
omit the details here, but it is worth noting that all this can be
handled nicely using a trick invented by the algebraic topologists
\cite{May}.  The idea is that there is a category Op whose objects are
opetopes.  The morphisms in this category describe how one opetope is
included as a specified face of another.  An opetopic set may then be
defined as a contravariant functor $S \maps {\rm Op} \to {\rm Set}$.
Such a functor assigns to each opetope $t$ a set $S(t)$ of cells of
that shape, or `$t$-cells'.  Moreover, if $f \maps s \to t$ is a
morphism in Op describing how $s$ is a particular face of $t$, the face
map $S(f) \maps S(t) \to S(s)$ describes how each $t$-cell of $S$ has a
given $s$-cell as this particular face.

For the definition of a weak $n$-category we need some terminology
concerning opetopic sets.  If $j \ge 1$, we may schematically represent
a $j$-dimensional cell $x$ in an opetopic set as follows:
\[ 
\begin{diagram}[(a_1,\dots,a_k)]
\node{(a_1,\dots,a_k)}\arrow{e,t}{x} \node{a'}
\end{diagram} 
\]
Here $a_1,\dots,a_k$ are the infaces of $x$ and $a'$ is the outface of
$x$; all these are cells of one lower dimension.  
A configuration just like this, satisfying all the incidence
relations satisfied by the boundary of a cell, but with $x$ itself
missing:
\[
\begin{diagram}[(a_1,\dots,a_k)]
\node{(a_1,\dots,a_k)}\arrow{e,t}{?} \node{a'}
\end{diagram}
\] 
is called a `frame'.  A `niche' is like a frame with the outface missing:
\[ 
\begin{diagram}[(a_1,\dots,a_k)]
\node{(a_1,\dots,a_k)}\arrow{e,t}{?} \node{?}
\end{diagram} 
\]
Similarly, a `punctured niche' is like a frame with the outface
and one inface missing:
\[ 
\begin{diagram}[(a_{1},a_{i+1},\dots,a_{k})]
\node{(a_{1},\dots,a_{i-1},?,a_{i+1},\dots,a_{k})}\arrow{e,t}{?} \node{?}
\end{diagram} 
\]
If one of these configurations (frame, niche,
or punctured niche) can be extended to an actual cell, the
cell is called an `occupant' of the configuration.  Occupants
of the same frame are called `frame-competitors', while occupants of
the same niche are called `niche-competitors'.  

\subsection{Universality} 

The only thing we need now to define the notion of weak $n$-category is
the concept of a `universal' occupant of a niche.  This is also the
subtlest aspect of the whole theory.  We explain it briefly here,
but it seems that the only way to really understand it is to carefully
work through examples.  

Before confronting the precise definition of universality, it is
important to note that the main role of universality is to define the
notion of `composite':

\begin{defn}\et  Given a universal occupant $u$ of a $j$-dimensional niche:
\[
\begin{diagram}[(a_1,\dots,a_k)]
\node{(a_1,\dots,a_k)}\arrow{e,t}{u} \node{b}
\end{diagram}
\]
we call $b$ a `composite' of $(a_{1},\dots,a_{k})$.
\end{defn}

It is also important to keep in mind the role played by cells of
different dimensions.  In our framework an $n$-category usually has
cells of arbitrarily high dimension.  For $j \le n$ the $j$-dimensional
cells play the role of $j$-morphisms, while for $j > n$ they play the
role of `equations', `equations between equations', and so on.  The
definition of universality depends on $n$ in a way that has the
following effects.  For $j \le n$ there may be many universal occupants
of a given $j$-dimensional niche, which is why we speak of `a' composite
rather than `the' composite.  There is at most one occupant of any given
$(n+1)$-dimensional niche, which is automatically universal.  Thus
composites of $n$-cells are unique, and we may think of the universal
occupant of an $(n+1)$-dimensional niche as an equation saying that the
composite of the infaces equals the outface.  For $j > n + 1$ there is
exactly one occupant of each $j$-dimensional frame, indicating that the
composite of the equations corresponding to the infaces equals the
equation corresponding to the outface.

The basic idea of universality is that a $j$-dimensional niche-occupant
is universal if all of its niche-competitors factor through it uniquely,
{\it up to equivalence}.  For $j \ge n+1$ this simply amounts to saying
that each niche has a unique occupant, while for $j = n$ it means that
each niche has an occupant through which all of its niche-competitors
factor uniquely.  In general, we require that composition with a
universal niche-occupant set up a `balanced punctured niche' of one
higher dimension.  Heuristically, one should think of a balanced
punctured niche as defining an equivalence between occupants of its
outface and occupants of its missing outface.

\begin{defn}\et A $j$-dimensional niche-occupant:
\[
\begin{diagram}[(c_1,\dots,c_k)]
\node{(c_1,\dots,c_k)}  \arrow{e,t}{u} \node{d}
\end{diagram}
\]
is said to be `universal' if and only if $j>n$ and $u$ is the only
occupant of its niche, 
or $j\le n$ and for any frame-competitor $d'$ of $d$, the
$(j+1)$-dimensional punctured niche
\[
\begin{diagram}[((c_{1},\dots,c_{k}) \to d,\; d \to d')]
\node{((c_{1},\dots,c_{k}) \mapright{u} d,\; d \mapright{?} d')} 
\arrow{s,r}{?} \\
\node{(c_{1},\dots,c_{k}) \mapright{?} d'}  
\end{diagram}
\]
and its mirror-image version
\[
\begin{diagram}[((c_{1},\dots,c_{k}) \to d,\; d \to d')]
\node{(d \mapright{?} d',\;(c_{1},\dots,c_{k}) \mapright{u} d)} 
\arrow{s,r}{?} \\
\node{(c_{1},\dots,c_{k}) \mapright{?} d'}  
\end{diagram}
\]
are balanced.  
\end{defn}

Finally we must define the concept of `balanced punctured niche'.  The
reader may note that the first numbered condition in the following
definition generalizes the concept of an `essentially surjective'
functor, while the second generalizes the concept of a `fully faithful'
functor.

\begin{defn} \et An $m$-dimensional punctured niche:
\[ 
\begin{diagram}[(a_{1},a_{j+1},\dots,a_{k})]
\node{(a_{1},\dots,a_{i-1},?,a_{i+1},\dots,a_{k})}\arrow{e,t}{?} \node{?}
\end{diagram}
\]
is said to be `balanced' if and only if $m>n+1$ or:
\begin{enumerate}
\item any extension
\[ 
\begin{diagram}[(a_{1},a_{i+1},\dots,a_{m})]
\node{(a_{1},\dots,a_{i-1},?,a_{i+1},\dots,a_{k})}\arrow{e,t}{?} \node{b}
\end{diagram} 
\]
extends further to:
\[ 
\begin{diagram}[(a_{1},a_{i+1},\dots,a_{k})]
\node{(a_{1},\dots,a_{i-1},a_i,a_{i+1},\dots,a_{k})}\arrow{e,t}{u} \node{b}
\end{diagram} 
\]
with $u$ universal in its niche, and 
\item  for any occupant
\[ 
\begin{diagram}[(a_{1},a_{i+1},\dots,a_{k})]
\node{(a_{1},\dots,a_{i-1},a_i,a_{i+1},\dots,a_{k})}\arrow{e,t}{u} \node{b}
\end{diagram} 
\]
universal in its niche, and frame-competitor $a'_{i}$ of $a_{i}$, the
$(m+1)$-dimensional punctured niche
\[
\begin{diagram}
[((c_{1},\dots,c_{m}) \to d,\; d \to d')]
\node{(a'_{i} \mapright{?} a_{i}, \;
(a_{1},\dots,a_{i-1},a_{i},a_{i+1},\dots,a_{k}) \mapright{u} b)}  
\arrow{s,r}{?}  \\
\node{(a_{1},\dots,a_{i-1},a'_{i},a_{i+1},\dots, a_{k}) \mapright{?} b}
\end{diagram}
\]
and its mirror-image version
\[
\begin{diagram}
[((c_{1},\dots,c_{m}) \to d,\; d \to d')]
\node{((a_{1},\dots,a_{i-1},a_{i},a_{i+1},\dots,a_{k}) \mapright{u} b, \;
a'_{i} \mapright{?} a_{i})} 
\arrow{s,r}{?}  \\
\node{(a_{1},\dots,a_{i-1},a'_{i},a_{i+1},\dots, a_{k}) \mapright{?} b}
\end{diagram}
\]
are balanced.  
\end{enumerate}
\end{defn}

Note that while the definitions of `balanced' and `universal' call upon
each other recursively, there is no bad circularity.  Using these
definitions, it is easy to define a weak $n$-category.  While the
definition below does not explicitly depend on $n$, it depends on $n$
through the definition of `universal' niche-occupant.

\begin{defn} \et A `weak $n$-category' is an opetopic
set such that 1) every niche has a universal occupant, and 2) composites
of universal cells are universal.  \end{defn}

\section{Conclusions}

The above definition of weak $n$-category is really a beginning, rather
than an end.  We turn the reader to the papers by Dolan and the author
\cite{BD2} and the forthcoming work of Hermida, Makkai and Power 
for more.  The weak $(n+1)$-category of weak $n$-categories is
beginning to be understood; the generalizations of $n$-cateories where
we replace $I$ by an arbitrary operad have also turned out to be very
interesting.   However, before the really interesting applications of
$n$-category theory can be worked out, there is still much basic work 
to be done.

In particular, it is important to compare various different definitions
of weak $n$-category, so that the subject does not fragment.  As one
might expect, the question of when two two definitions of weak
$n$-category are `equivalent' is rather subtle.  This question seems to
have first been seriously pondered by Grothendieck \cite{Gro}, who
proposed the following solution.  Suppose that for all $n$ we have two
different definitions of weak $n$-category, say `$n$-${\rm category}_1$'
and `$n$-${\rm category}_2$'.  Then we should try to construct the
$(n+1)$-${\rm category}_1$ of all $n$-${\rm categories}_1$ and the
$(n+1)$-${\rm category}_1$ of all $n$-${\rm categories}_2$ and see if
these are equivalent as objects of the $(n+2)$-${\rm category}_1$ of all
$(n+1)$-${\rm categories}_1$.  If so, we may say the two definitions are
equivalent as seen from the viewpoint of the first definition.

There are some touchy points here worth mentioning.  First, there is
considerable freedom of choice involved in constructing the two
$(n+1)$-${\rm categories}_1$ in question; one should do it in a
`reasonable' way, but this is not necessarily easy.  Secondly, there is
no guarantee that we might not get a different answer for the question
if we reversed the roles of the two definitions.  Nonetheless, it should
be interesting to compare different definitions of weak $n$-category in
this way.

A second solution is suggested by homotopy theory, which again comes to
the rescue.  Many different approaches to homotopy theory are in use,
and though superficially very different, there is a well-understood
sense in which they are fundamentally the same.  Different approaches
use objects from different categories to represent topological spaces,
or more precisely, the homotopy-invariant information in topological
spaces, called their `homotopy types'.  These categories are not
equivalent, but each one is equipped with a class of morphisms playing
the role of homotopy equivalences.  Given a category $C$ equipped with a
specified class of morphisms called `equivalences', under mild
assumptions one can adjoin inverses for these morphisms, and obtain a
category called the the `homotopy category' of $C$.  Two categories with
specified equivalences may be considered the same for the purposes of
homotopy theory if their homotopy categories are equivalent in the usual
sense of category theory.  Homotopy theorists have proved that all the
popular approaches to homotopy theory are the same in this sense
\cite{Baues}.

The same strategy should be useful in $n$-category theory.  Any
definition of weak $n$-category should come along with a definition of
an `$n$-functor' for which there is a category with weak $n$-categories
as objects and $n$-functors as morphisms, and there should be a
specified class of $n$-functors called `equivalences'.  This allows the
construction a homotopy category of $n$-categories.  Then, for two
definitions of weak $n$-category to be considered equivalent, we
require that their homotopy categories be equivalent.

Dolan and the author have constructed the homotopy category of their
$n$-categories, and Simpson \cite{Simpson} has constructed the homotopy
category of Tamsamani's $n$-categories.  Now we need machinery to check
whether these homotopy categories, and those corresponding to other
definitions, are equivalent.  Once these preliminary chores are
completed, there should be many exciting things we can do with
$n$-categories.

\end{document}